\documentclass[a4paper,oneside]{amsart}
\usepackage[english]{babel}
\usepackage[applemac]{inputenc} 
\usepackage[T1]{fontenc}
\usepackage{mathpazo}

\usepackage[pagewise,displaymath]{lineno}

\usepackage{titlesec} 
\titleformat{\section}{\center\Large}{\Large\thesection.}{.5em}{}
\titleformat{\subsection}{\center\large}{\large \thesubsection}{.5em}{}
\titleformat{\subsubsection}{\center\large}{\large \thesubsubsection}{.5em}{}

\usepackage[leqno]{amsmath}
\usepackage{amssymb}
\usepackage{stmaryrd}             
\usepackage{amsthm}        
\usepackage{dsfont}       
\usepackage{upref}                   
\usepackage{amsfonts}             
\usepackage{exscale}
\usepackage{textcomp}
\usepackage{paralist}       
\usepackage{ragged2e}     
\usepackage{varioref}         
\usepackage{calc}
\titlelabel{\thesection. \enspace}

\newcommand\R{{\mathds{R}}}        

\newcommand\abs[1]{\lvert#1\rvert} 
\newcommand\norm[1]{\lVert#1\rVert} 
\newtheoremstyle{prefstyle}
  {14pt}
  {14pt}
  {\normalfont}         
  {\parindent}
  {\scshape\bfseries}
  {.}
  {7pt}
  {}
  \theoremstyle{prefstyle}
 \newtheorem{defn}{Definition}[section]

  \newtheorem{rem}[defn]{Remark}
 
\newtheoremstyle{mystylethm}
  {14pt}
  {14pt}
  {\itshape}
  {\parindent}
  {\bfseries\scshape}
  {.}
  {7pt}
  { }

\theoremstyle{mystylethm}       
\newtheorem{thm}{Theorem}[section]
 \newtheorem{cor}[thm]{Corollary}
  \newtheorem{lem}[defn]{Lemma}
\newtheoremstyle{Appstylethm}
  {14pt}
  {14pt}
  {\itshape}
  {\parindent}
  {\bfseries\scshape}
  {.}
  {7pt}
  { }

\newtheorem{appthm}{Theorem A.}

\newenvironment{myproof}
                           {{\fontseries{m}\scshape\bfseries
                           \selectfont Proof.\;}}{\qed\\[7pt]}

\newenvironment{proofthm}[1] 
                           {{\fontseries{m}\scshape\bfseries
                           \selectfont Proof of Theorem #1.\;}}
                           {\qed\\[7pt]}
                           
                           {{\fontseries{m}\scshape\bfseries
                           \selectfont Proof of Lemma #1.\;}}
                           {\qed\\[7pt]}

\newenvironment{proofcor}[1] 
                           {{\fontseries{m}\scshape\bfseries
                           \selectfont Proof of Corollary #1.\;}}
                           {\qed\\[7pt]}

\def\XXint#1#2#3{{\setbox0=\hbox{$#1{#2#3}{\int}$ }
\vcenter{ \hbox{$#2#3$} }\kern-.45\wd0}}

\fontsize{14}{14}\selectfont

\begin{document}

\thispagestyle{empty}
\title[New weighted Hardy's inequalities with application to
nonexistence of global solutions]{New weighted Hardy's inequalities with application to
nonexistence of global solutions}

\author{By Daniel Hauer, Abdelaziz Rhandi} 

\address{Institut f\"ur Angewandte Analysis, Universit\"at Ulm, 89069 Ulm,
Germany. \phantom{,}\newline \mbox{}\hspace{9pt} 
Universit\'e de Lorraine (campus de Metz) et CNRS, Laboratoire de
Math\'ematiques et Applications de Metz, UMR 7122, B\^at. A, Ile du Saulcy,
57045 Metz Cedex 1, France.}
\email{daniel.hauer@uni-ulm.de, hauer@math.cnrs.fr.}

\address{Dipartimento di Matematica,
Universit\`a degli Studi di Salerno, via Ponte don Melillo, 84084 Fisciano (Sa),
Italy.}
\email{arhandi@unisa.it.}
\date{\today} 



\begin{abstract}
In this article, we prove that the following weighted Hardy inequality
\begin{equation} \label{hardy-inequality:1}
\big(\tfrac{\abs{d-p}}{p}\big)^{p}\,\int\, \tfrac{\abs{u}^{p}}{\abs{x}^{p}}\; d\mu
\le \int\,\abs{\nabla u}^{p}\;d\mu + \big(\tfrac{\abs{d-p}}{p}\big)^{p-1}\,
\textrm{sgn}(d-p)\,  \int\, \abs{u}^{p}\,
\tfrac{(x^{t}Ax)^{p/2}}{\abs{x}^{p}}\; d\mu
\end{equation}
holds for all $u$ in the weighted Sobolev space $W^{1,p}_{\mu}$ with best constant,
where $A\in \R^{d\times d}$ is a positive-definite symmetric matrix, $d\ge1$, $1<p<+\infty$,
and $\mu$ denotes a Borel measure on $\R^{d}$ given by 
\begin{equation}\label{m:mu}
d\mu=\rho(x)\,dx\qquad\text{with density}\qquad \rho(x)=c\cdot\exp(-\tfrac{1}{p}(x^{t}Ax)^{p/2}),
\quad (c>0).
\end{equation}
Here the integral is taken over $]0,+\infty[$
if $d=1$ and over $\R^{d}$ if $d\ge2$. If $p>d$, then we can
deduce from inequality~\eqref{hardy-inequality:1}, that there is a Poincar\'e inequality on
$W^{1,p}_{\mu}$. The proof of inequality~\eqref{hardy-inequality:1} is based on 
the method of vector fields firstly introduced by Mitidieri~\cite{MR1769903}.
By the same method, we prove for the same weights and for $1<p<+\infty$, $d\ge1$,
weighted Caffarelli-Kohn-Nirenberg inequalities. As an application of inequality
\eqref{hardy-inequality:1} we prove a nonexistence result for a p-Kolmogorov
parabolic equation.
\end{abstract}
\maketitle
%
\setcounter{section}{0}
\section{{\bf Introduction}}
\mbox{}\indent In the cases, when $A\equiv 0$, and $c=1$, the
Borel measure $\mu$ defined in~\eqref{m:mu} reduces to the Lebesgue
measure on $\R^{d}$ and so inequality
\eqref{hardy-inequality:1} becomes the well-known Hardy inequality 
\begin{equation}\label{hardy:classic}
\Big(\tfrac{\abs{d-p}}{p}\Big)^{p}\,\int \tfrac{\abs{u}^{p}}{\abs{x}^{p}}\,dx\le 
\int \abs{\nabla u}^{p}\,dx\qquad (u\in W^{1,p})\mathpunct{.}
\end{equation}
Inequality~\eqref{hardy:classic} was first stated in dimension $d=1$
in~\cite{MR1544414} by Hardy in 1920. Various generalization of Hardy's inequality
since have been found with application to various branches of mathematics,
see for instance \cite[p. 175f]{MR1190927} by Mitrinovi\'c, Pe\v{c}ari\'c, and Fink.
Among them, we want shortly summarize the relation between the
Hardy inequality with its optimal constant and the existence and nonexistence
theory of nonnegative solutions of parabolic equations containing a critical
potential. In 1984 Baras and J. Goldstein proved in~\cite{MR742415}
the following important result. 
\setcounter{appthm}{0}
\begin{appthm}{{\bf (Baras-Goldstein)}}\label{thm:A}
Let $\Omega=]0,+\infty[$ if $d=1$ and $\Omega=\R^{d}$ if $d\ge2$.
\begin{enumerate}
\item[(i)] If $\lambda\le \big(\tfrac{d-2}{2}\big)^{2}$, then for the linear
parabolic equation
\begin{equation}\label{KV:2}
\tfrac{\partial u}{\partial t}-\Delta u= \tfrac{\lambda}{\abs{x}^{2}}u\qquad x\in
\Omega\mathpunct{,}\;t>0\mathpunct{,}
\end{equation}
admits a nonnegative nontrivial solution.
\item[(ii)] If $\lambda>\big(\tfrac{d-2}{2}\big)^{2}$, then
equation~\eqref{KV:2} admits no nonnegative nontrivial solution.
\end{enumerate}
\end{appthm}
Obviously, the phenomenon of existence and nonexistence is caused by
the singular potential $\lambda \abs{x}^{-2}$, which is controlled by
Hardy's inequality \eqref{hardy:classic} together with its optimal constant.
Improvements of inequality~\eqref{hardy:classic} with or without reminder term
of inequality~\eqref{hardy:classic} have been found, 
see for instance~\cite[1997]{MR1605678} by Brezis and V\'azquez,
\cite[1998]{MR1616905} by Garc\'ia Azorero and Peral Alonso,
\cite[2000]{MR1760280} by V\' azquez and Zuazua,
\cite[2006]{Arendt-Goldstein-2006} by W. Arendt, G. R. Goldstein, and J. A.
Goldstein, as well as~\cite[2011]{Goldstein-Rhandi-weighted-hardy-11} by G. R.
Goldstein, J. A. Goldstein and Rhandi, and served them in their studies about
existence, nonexistence, and the qualitative behavior of solutions. In this way, for
instance, the authors in~\cite{Goldstein-Rhandi-weighted-hardy-11},
have established first inequality \eqref{hardy-inequality:1} when $d\ge3$,
$p=2$, $A\in \R^{d\times d}$ is a
positive-definite symmetric matrix, and $c$ is such that $\int\,\rho\,dx=1$. Then, they proved the following result:
\begin{appthm}{{\bf (Goldstein-Goldstein-Rhandi)}}\label{thm:B}
\begin{enumerate}
\item[(i)] If $\lambda\le \big(\tfrac{d-2}{2}\big)^{2}$, then for every nonnegative
initial value $u_{0}\in L^{2}_{\mu}$, equation
\begin{equation}\label{KV:3}
\tfrac{\partial u}{\partial t}-\Delta u= \langle Ax,\nabla u\rangle 
+ \tfrac{\lambda u}{\abs{x}^{2}}\quad x\in \R^{d},\;t>0
\end{equation}
admits a nonnegative global solution with exponential growth.
\item[(ii)] If $\lambda>\big(\tfrac{d-2}{2}\big)^{2}$, then for every
nonnegative initial value $u_{0}\in L^{2}_{\mu}\setminus\{0\}$ there is no
nonnegative global solution of equation~\eqref{KV:3} with exponential growth.
\end{enumerate}
\end{appthm}
In order to prove the nonexistence of global solutions of
equation~\eqref{KV:3}, the authors in
\cite{Goldstein-Rhandi-weighted-hardy-11} employ an
approach in~\cite{MR1733904} due to Cabr\'e and Martel. Comparing
Theorem~\ref{thm:B} with Theorem~\ref{thm:A}, one sees that the
unbounded drift term $Ax$ appearing in the symmetric Ornstein-Uhlenbeck operator
$L=-\Delta+\langle Ax,\nabla \cdot\rangle$, has a strong influence
on the qualitative behavior on the corresponding
solutions of equation~\eqref{KV:3}.\\[7pt]
\mbox{}\indent In analogs with the two theorems above, we first state and prove in
Section 2 of this article the weighted Hardy inequality~\eqref{hardy-inequality:1}. If
$p>d$, then we deduce from inequality~\eqref{hardy-inequality:1} a weighted
Poincar\'e inequality. The proof of inequality~\eqref{hardy-inequality:1} is here
based on the so-called method of vector field first introduced by Mitidieri in
\cite{MR1769903}. Based on the same method, we state weighted
Caffarelli-Kohn-Nirenberg inequalities in Section 3. In the last section, we apply
our Hardy inequality and use the optimality of the constant
$\big(\tfrac{\abs{d-p}}{p}\big)^{p}$
to prove an existence and nonexistence result in dimension $d=1$ of
equation~\eqref{KV:3}, where we replace the symmetric Ornstein-Uhlenbeck
operator $L$ in~\eqref{KV:3} by the nonlinear p-Kolmogorov operator
\[
K_{p}u=-\Delta_{p}+\langle\abs{\nabla\cdot}^{p-2}\nabla\cdot,\tfrac{\nabla\rho}{\rho}
\rangle\mathpunct{.}
\]
The p-Kolmogorov operator was first introduced by G. R. Goldstein, J. A.
Goldstein and Rhandi in~\cite{MR2746424}. We will use an idea of
Cabr\'e and Martel in~\cite{MR1733904}.\\[7pt] 
\mbox{}\indent Throughout this article,
we employ the following notation: For $1\le p\le +\infty$ and if $\mu$ is
the Borel measure defined by \eqref{m:mu}, then we denote by
$L^{p}_{\mu}:=L^{p}(d\mu)$ the weighted Lebesgue space on either $]0,+\infty[$
if $d=1$ or on $\R^{d}$ if $d\ge2$, and by $L^{p}_{\mu,loc}$ the set of locally 
$L^{p}_{\mu}$-integrable function on either $]0,+\infty[$ if $d=1$ or on $\R^{d}$
if $d\ge2$. Further, we denote by $W^{1,p}_{\mu}$ the first weighted
Sobolev space, i.e., the set of all $u\in L^{p}_{\mu}$ having the distributional
partial derivative $\partial_{x_{1}}u, \dots, \partial_{x_{d}}u\in L^{p}_{\mu}$.
Similarly, an element $u\in W^{1,p}_{\mu,loc}$ if and only if $u\in L^{p}_{\mu,loc}$
and $\partial_{x_{1}}u, \dots, \partial_{x_{d}}u\in L^{p}_{\mu,loc}$. In this context,
it is worth mentioning, that for the measure $\mu$ defined by~\eqref{m:mu},
the set $\textrm{C}_{c}^{\infty}$ of infinitely differentiable function having compact
support in either $]0,+\infty[$ if $d=1$ or on $\R^{d}$ if $d\ge2$ lies dense in
$W^{1,p}_{\mu}$, see for instance ~\cite{Tolle:2011yu} by T\"olle.
\newpage
%
%
\section{{\bf Main results}}
In this section, we study for $1<p<+\infty$ the following weighted Hardy inequality with optimal constant.

\begin{thm}\label{thm-hardy-inequality}
Let $A\in \R^{d\times d}$ be a symmetric positive
semi-definite matrix. Then, for all $u\in W^{1,p}_{\mu}$,
\begin{equation}\tag{\ref{hardy-inequality:1}}
\Big(\tfrac{\abs{d-p}}{p}\Big)^{p}\,\int\, \tfrac{\abs{u}^{p}}{\abs{x}^{p}}\; d\mu
\le \int\,\abs{\nabla u}^{p}\;d\mu + \Big(\tfrac{\abs{d-p}}{p}\Big)^{p-1}\,
\textrm{sgn}(d-p)\,  \int\, \abs{u}^{p}\,
\tfrac{(x^{t}Ax)^{\frac{p}{2}}}{\abs{x}^{p}}\; d\mu\mathpunct{.}
\end{equation}
Moreover, if either $A\equiv0$ or $A$ is positive definite, then the constant
$C(d,p)=\Big(\tfrac{\abs{d-p}}{p}\Big)^{p}$ is optimal (including the
case $p=d$ if $d\ge2$). Here the integral is taken over $]0,+\infty[$ if $d=1$ and
over $\R^{d}$ if $d\ge2$.
\end{thm}

\begin{rem}
Since for any symmetric positive semi-definite matrix 
$A=\in \R^{d\times d}$, we have that $(x^{t}Ax)\le \abs{x}^{2}\,\abs{A}$ for all $x\in \R^{d}$, one easily
sees that the second term on the right hand-side in~\eqref{hardy-inequality:1}
satisfies
\[
\textrm{sgn}(d-p)\,  \int\, \abs{u}^{p}\,
\tfrac{(x^{t}Ax)^{\frac{p}{2}}}{\abs{x}^{p}}\; d\mu
\le \abs{A}^{\frac{p}{2}}\,\int \abs{u}^{p}\,d\mu\qquad\text{for all $u\in L^{p}_{\mu}$.}
\]
Thus Hardy's inequality~\eqref{hardy-inequality:1} implies that
$W^{1,p}_{\mu}\hookrightarrow L^{p}(\tfrac{C(d,p)}{\abs{x}^{p}}\,d\mu)$ by a continuous injection, provided $p\neq d$. Furthermore, we can deduce from inequality~\eqref{hardy-inequality:1} the following
Poincar\'e inequality.
\end{rem}

\begin{cor}\label{cor:poincare-inequality}
If the matrix $A$ is positive-definite, then for $p>d$,
\begin{equation}\label{poincare-inequality}
\Big(\tfrac{p-d}{p}\Big)^{p-1}\,\lambda^{p/2}(A)\, \int\, \abs{u}^{p}\,d\mu\le \int\,\abs{\nabla u}^{p}\;d\mu\qquad\text{for all $u\in W^{1,p}_{\mu}$,}
\end{equation}
where $\lambda(A)>0$ denotes the lowest eigenvalue of $A$.
\end{cor}

\begin{proofcor}{\ref{cor:poincare-inequality}}
If the symmetric matrix $A\in \R^{d}$ is positive definite, then the lowest eigenvalue
$\lambda(A)$ of $A$ is strictly positive and $x^{t}Ax\ge \lambda(A)\,\abs{x}^{2}$ for all $x\in \R^{d}$,
whence we can deduce inequality~\eqref{poincare-inequality} from inequality~\eqref{hardy-inequality:1}
provided $p>d$.
\end{proofcor}

\begin{proofthm}{\ref{thm-hardy-inequality}}
It is not hard to see that in the case $p=d$ for $d\ge2$,
inequality~\eqref{hardy-inequality:1} holds true with optimal constant $C(d,d)=0$.
It is left to show that inequality~\eqref{hardy-inequality:1} holds true when
$p\neq d$ and $d\ge1$. To do so, we follow an approach introduced by Mitidieri
in~\cite{MR1769903}. We take $\varepsilon>0$, $\lambda\ge0$, which will be
chosen later, and set
\[
F(x)= \lambda\,\textrm{sgn}(d-p)\, \frac{x}{\abs{x}^{p}+\varepsilon}\,\rho(x)
\]
for every $x\in\R^{d}$ if $d\ge2$ and for all $x\in ]0,+\infty[$ if $d\ge1$, where
$\rho$ is defined in~\eqref{m:mu}. Since for every $i=1,\dots, d$,
\[
\frac{\partial}{\partial x_{i}}\frac{x_{i}}{\abs{x}^{p}+\varepsilon}
=\frac{1}{\abs{x}^{p}+\varepsilon} -\frac{p\, x_{i}^{2}\,\abs{x}^{p-2}}{(\abs{x}^{p}+\varepsilon)^{2}}
\]
and since for $f(x):=x^{t}Ax$, $(x\in \R^{d})$, $f'(x)=x^{t}(A+A^{t})=2x^{t}A$,
we have that
\[
\frac{\partial F_{i}}{\partial x_{i}} = \lambda\,\textrm{sgn}(d-p)\, 
 \left[\frac{1}{\abs{x}^{p}+\varepsilon}
-\frac{p\, x_{i}^{2}\,\abs{x}^{p-2}}{(\abs{x}^{p}+\varepsilon)^{2}} - 
\frac{x_{i}(x^{t}A)_{i}}{\abs{x}^{p}+\varepsilon}  (x^{t}Ax)^{\frac{p}{2}-1}  \right]\,
\rho(x)
\]
and so
\[
\textrm{div}\big( F(x) \big) = \lambda\, \,\textrm{sgn}(d-p)\,
 \left[\frac{d}{\abs{x}^{p}+\varepsilon}
-\frac{p\,\abs{x}^{p}}{(\abs{x}^{p}+\varepsilon)^{2}} - 
\frac{ (x^{t}Ax)^{\frac{p}{2}}}{\abs{x}^{p}+\varepsilon}  \right]
 \, \rho(x)\mathpunct{.}
\]
Since the set of infinitely differentiable functions with compact support lies dense
in $W^{1,p}_{\mu}$, it is sufficient to show that
inequality~\eqref{hardy-inequality:1} holds for $u\in \textrm{C}_{c}^{\infty}$. Fix
$u\in \textrm{C}_{c}^{\infty}$. Then, by an integration by parts and 
Young's inequality,
\allowdisplaybreaks 
\begin{align*}
\int\,\abs{u}^{p} \,\lambda\,\textrm{sgn}(d-p)\, \left[\tfrac{d}{\abs{x}^{p}+\varepsilon}
-\tfrac{p\,\abs{x}^{p}}{(\abs{x}^{p}+\varepsilon)^{2}} - 
\tfrac{ (x^{t}Ax)^{\frac{p}{2}}}{\abs{x}^{p}+\varepsilon}\right] \,d\mu
& =  \int\,\abs{u}^{p} \textrm{div}\big( F \big)\; dx\\
& = (-p) \int\,\abs{u}^{p-1}\;\textrm{sgn}(u)\,\langle\nabla u,F\rangle\; dx\\
& \le \int\,\abs{\nabla u}^{p}\;d\mu + (p-1)\,\lambda^{p\prime} 
\,\int\,\tfrac{\abs{u}^{p}}{\abs{x}^{p}}\; d\mu\mathpunct{.}
\end{align*}
And hence
\begin{align*}
&\int\,\abs{u}^{p} \,\lambda\,\textrm{sgn}(d-p)\, \left[\tfrac{d}{\abs{x}^{p}+\varepsilon}
-\tfrac{p\,\abs{x}^{p}}{(\abs{x}^{p}+\varepsilon)^{2}}\right]\,d\mu - 
(p-1)\,\lambda^{p\prime} \,\int\,\tfrac{\abs{u}^{p}}{\abs{x}^{p}} \,d\mu\\
& \qquad
\le \int\,\abs{\nabla u}^{p}\;d\mu +\lambda\,\textrm{sgn}(d-p)\,
\int\,\abs{u}^{p} \tfrac{ (x^{t}Ax)^{\frac{p}{2}}}{\abs{x}^{p}+\varepsilon}d\mu
\mathpunct{.}
\end{align*}
Letting $\varepsilon\to0+$, applying Lebesgue's dominated convergence
theorem and Fatou's lemma, we obtain that
\[
\int\, \tfrac{\abs{u}^{p}}{\abs{x}^{p}} \big[\lambda\, \abs{d-p} - (p-1)\,\lambda^{p\prime}\big]\; d\mu
\le \int\,\abs{\nabla u}^{p}\;d\mu + \lambda\,\textrm{sgn}(d-p)\,  
\int \abs{u}^{p}\,\tfrac{(x^{t}Ax)^{\frac{p}{2}}}{\abs{x}^{p}}\; d\mu
\mathpunct{.}
\]
Now, we choose $\lambda=\Big(\tfrac{\abs{d-p}}{p}\Big)^{p-1}$ in the last inequality,
which is, in fact, the maximum of the function $\lambda\mapsto \big[\lambda\,\abs{d-p} - (p-1)\,
\lambda^{p\prime}\big]$ on the half line $[0,+\infty[$ and achieve to
inequality~\eqref{hardy-inequality:1}.\\[7pt]
%
\mbox{}\indent Next, we show the optimality of the constant
$\big(\tfrac{\abs{d-p}}{p}\big)^{p}$ when $d\neq p$ and when the matrix $A$
is positive definite. To do so, let $\lambda>\big(\tfrac{\abs{d-p}}{p}\big)^{p}$, and for
$\gamma$ such that
\begin{equation}\label{choosen-gamma}
1-\tfrac{d}{p} <\gamma < 0 \quad\text{ if $p<d$}\qquad\text{and}\qquad
1-\tfrac{d}{p} < \gamma < 1 \quad\text{ if $p>d$,}
\end{equation}
set $\varphi(x)=\abs{x}^{\gamma}$. Then, $\nabla\varphi(x)=\gamma\, 
\abs{x}^{\gamma-2}x$ for all $x\neq 0$, and since $x^{t}Ax\le \abs{A}\,\abs{x}^{2}$
for all $x\in \R^{d}$, we have that
\begin{equation}\label{form:1}
\begin{split}
&\;\int\,\Big[\abs{\nabla\varphi}^{p} + \lambda^{1/p\prime}\,
\textrm{sgn}(d-p)\,  \abs{\varphi}^{p}\,\frac{(x^{t}Ax)^{\frac{p}{2}}}{\abs{x}^{p}}\;
-\frac{\lambda}{\abs{x}^{p}}\abs{\varphi}^{p} \Big]\; d\mu\\
&\hspace{2cm} \le 
[\abs{\gamma}^{p}-\lambda]\, \int\, \abs{x}^{p(\gamma-1)}\,
d\mu + \abs{A}^{p/2}\,\lambda^{\frac{1}{p\prime}} \int\,\abs{x}^{\gamma\, p}\,
 d\mu\mathpunct{.} 
\end{split}
\end{equation}
By the assumption, there are $\alpha_{1}, \alpha_{2}>0$ such that
\[
\alpha_{1}\,\abs{x}^{2}\le x^{t}Ax \le \alpha_{2}\,\abs{x}^{2}\qquad\text{for all $x\in \R^{d}$.}
\]
Hence for every $\beta\in \R$,
\begin{equation}\label{form:2}
\int\, \abs{x}^{p\beta}\,d\mu
 \le  c\,\int_{0}^{+\infty} \abs{x}^{p\beta}\,
 e^{-\tfrac{\alpha_{1}^{p/2}\abs{x}^{p}}{p}}\,dx
 \quad
\text{and}\quad
\;\int\, \abs{x}^{p\beta}\,d\mu
 \ge  c\,\int_{0}^{+\infty}\abs{x}^{p\beta}\, e^{-\tfrac{\alpha_{2}^{p/2}\abs{x}^{p}}{p}}\,dx
 \mathpunct{.}
\end{equation}
For every $i=1, 2$, and every $\beta\in \R$,
\begin{equation}\label{form:4}
\begin{split}
\int\,\abs{x}^{p\beta}\,e^{-\tfrac{\alpha_{i}^{p/2}\abs{x}^{p}}{p}}\,dx
&= \,\sigma(S_{d-1})\,\int_{0}^{+\infty} r^{p\beta}\,
 e^{-\tfrac{\alpha_{i}^{p/2}\,r^{p}}{p}}\;r^{d-1}\,dr\\
&= \,\sigma(S_{d-1})\,p^{\frac{p\beta+d-p}{p}} \alpha_{i}^{-\frac{p\beta+d}{2}}
\int_{0}^{+\infty}\,t^{\beta+\frac{d}{p}-1}\,e^{-t}\;dt\mathpunct{,}
\end{split}
\end{equation}
where $\sigma(S_{d-1})$ denotes the total surface measure of the
unite sphere $S_{d-1}:=\{x\in \R^{d}\,\big\lvert\;\abs{x}=1\}$ with respect to the
surface measure $\sigma$ on $S_{d-1}$. We note that
\[
\int_{0}^{+\infty}\,t^{\beta+\frac{d}{p}-1}\,e^{-t}\;dt=\Gamma(\beta+\tfrac{d}{p})
\qquad\text{ is finite for every $\beta>-\tfrac{d}{p}$}
\]
and in particular for $\beta=\gamma$ or $\beta=\gamma-1$ when we choose
$\gamma$ as in~\eqref{choosen-gamma}. Thus, $\varphi(x)=\abs{x}^{\gamma}$
belongs to $W^{1,p}_{\mu}$, and in view of \eqref{form:1}-\eqref{form:4},
we have that
\allowdisplaybreaks
\begin{align*}
&\; \frac{\displaystyle\int\,\textstyle\Big[\abs{\nabla\varphi}^{p} + \lambda^{\frac{1}{p\prime}}\,
\textrm{sgn}(d-p)\,  \abs{\varphi}^{p}\,\frac{(x^{t}Ax)^{p/2}}{\abs{x}^{p}}\;
-\frac{\lambda}{\abs{x}^{p}}\abs{\varphi}^{p} \Big]\; d\mu}{\displaystyle
\int\textstyle\abs{\varphi}^{p}\,d\mu}\\
&\hspace{0,5cm} \le \frac{[\abs{\gamma}^{p}-\lambda]\,\displaystyle\int\,\textstyle\,
 \abs{x}^{p(\gamma-1)}\,d\mu
  + \abs{A}^{p/2}\,\lambda^{\frac{1}{p\prime}}\,\displaystyle\int\,\textstyle\,\abs{x}^{\gamma\, p}\, 
  d\mu}{\displaystyle\int\,\textstyle\, 
  \abs{x}^{p\gamma}\,d\mu}\\
& \hspace{0,5cm} \le  \frac{[\abs{\gamma}^{p}-\lambda]\,
p^{\frac{p(\gamma-1)+d-p}{p}} \alpha_{1}^{-\frac{p(\gamma-1)+d}{2}}
\Gamma((\gamma-1)+\frac{d}{p}) +\abs{A}^{p/2}\,\lambda^{\frac{1}{p\prime}}\,
 \,p^{\frac{p\gamma+d-p}{p}} \alpha_{1}^{-\frac{p\gamma+d}{2}}
\Gamma(\gamma+\frac{d}{p})}
{\,p^{\frac{p\gamma+d-p}{p}} \alpha_{2}^{-\frac{p\gamma+d}{2}}
\Gamma(\gamma+\frac{d}{p})}\\
& \hspace{0,5cm} =  \frac{[\abs{\gamma}^{p}-\lambda]\,
p^{\frac{p(\gamma-1)+d-p}{p}} \alpha_{1}^{-\frac{p(\gamma-1)+d}{2}} 
+\abs{A}^{p/2}\,\lambda^{\frac{1}{p\prime}}\,\,p^{\frac{p\gamma+d-p}{p}}
 \alpha_{1}^{-\frac{p\gamma+d}{2}}\,[(\gamma-1)+\frac{d}{p}]}
 {\,p^{\frac{p\gamma+d-p}{p}} \alpha_{2}^{-\frac{p\gamma+d}{2}}\,
 [(\gamma-1)+\frac{d}{p}]}
\end{align*}
Since, 
\begin{align*}
&\;\lim_{\gamma\to(1-\frac{d}{p})+} p^{\frac{p(\gamma-1)+d-p}{p}} = p^{-1}\mathpunct{,}
\qquad \lim_{\gamma\to (1-\frac{d}{p})+} p^{\frac{p\gamma+d-p}{p}} = 1\mathpunct{,}
\qquad
\lim_{\gamma\to (1-\frac{d}{p})+} \alpha_{1}^{-\frac{p(\gamma-1)+d}{2}} = 1\mathpunct{,}\\
&\; \lim_{\gamma\to (1-\frac{d}{p})+} \alpha_{i}^{-\frac{p\gamma+d}{2}}
 = \alpha_{i}^{-\frac{p}{2}}\;(i=1,2)\mathpunct{,}\qquad
\lim_{\gamma\to (1-\frac{d}{p})+} \Big[(\gamma-1) +\tfrac{d}{p}\Big]=0\mathpunct{,}
\end{align*}
and since
\[
\lim_{\gamma\to (1-\frac{d}{p})+}(\abs{\gamma}^{p}-\lambda) 
= \Big(\tfrac{\abs{d-p}}{p}\Big)^{p}-\lambda<0\mathpunct{,}\qquad
\]
we have  that
\[
\lim_{\gamma\to (1-\frac{d}{p})+}\;\;\frac{[\abs{\gamma}^{p}-\lambda]\,
p^{\frac{p(\gamma-1)+d-p}{p}} \alpha_{1}^{-\frac{p(\gamma-1)+d}{2}} 
+\abs{A}^{p/2}\,\lambda^{\frac{1}{p\prime}}\,\,p^{\frac{p\gamma+d-p}{p}}
 \alpha_{1}^{-\frac{p\gamma+d}{2}}\,[(\gamma-1)+\frac{d}{p}]}
 {\,p^{\frac{p\gamma+d-p}{p}} \alpha_{2}^{-\frac{p\gamma+d}{2}}\,
 [(\gamma-1)+\frac{d}{p}]}=-\infty\mathpunct{.}
\]
Thus,
\[
\inf_{\varphi \in W^{1,p}_{\mu}: \norm{\varphi}_{L^{p}_{\mu}}>0}
\frac{\displaystyle\int\,\abs{\nabla\varphi}^{p}\,d\mu +  \lambda^{\frac{1}{p\prime}}\,
\textrm{sgn}(d-p)\,\displaystyle\int \,\abs{\varphi}^{p}\,
\frac{(x^{t}Ax)^{p/2}}{\abs{x}^{p}}\,d\mu\;
- \displaystyle\int\, \frac{\lambda}{\abs{x}^{p}}\abs{\varphi}^{p} \, d\mu}{\displaystyle
\int\,\abs{\varphi}^{p}\,d\mu}=-\infty
\]
and hence, for every $M>0$, there is a $\varphi\in W^{1,p}_{\mu}$ with
$\norm{\varphi}_{L^{p}_{\mu}}>0$ satisfying
\[
\int\,\abs{\nabla\varphi}^{p}\,d\mu +  \lambda^{\frac{1}{p\prime}}\,
\textrm{sgn}(d-p)\,\int \,\abs{\varphi}^{p}\,
\frac{(x^{t}Ax)^{p/2}}{\abs{x}^{p}}\,d\mu\;
- \int\, \frac{\lambda}{\abs{x}^{p}}\abs{\varphi}^{p} \, d\mu 
< (-M) \int\, \abs{\varphi}^{p}\; d\mu< 0\mathpunct{.}
\]
This shows that the constant $\big(\tfrac{\abs{d-p}}{p}\big)^{p}$ in
inequality~\eqref{hardy-inequality:1} is optimal.
\end{proofthm}
%
%
\section{{\bf Two weighted Caffarelli-Kohn-Nirenberg inequalities}}

\mbox{}\indent We follow again the Ansatz as outlined in~\cite{MR1769903} and we
prove for $1<p<+\infty$ the following two weighted inequalities with optimal constant. Both inequalities
reduce in the case $p=1$, $c=1$, $A\equiv 0$ to the famous Caffarelli-Kohn-Nirenberg
inequalities~\cite{MR768824}. Here and as above, the integral is taken over $]0,+\infty[$ if $d=1$ and over $\R^{d}$ if $d\ge2$. We denote by
$\textrm{C}_{c,0}^{\infty}$ either the set $\textrm{C}_{c}^{\infty}(]0,+\infty[)$ if $d=1$
or $\textrm{C}_{c}^{\infty}(\R^{d}\setminus\{0\})$ if $d\ge2$.

\begin{thm}\label{caffarelli:1}
Let $A\in \R^{d\times d}$ be a symmetric positive
semi-definite matrix, and let $a\in\R$. Then 
\[
\Big(\tfrac{\abs{d-p(a+1)}}{p}\Big)^{p}\,\int\, \tfrac{\abs{u}^{p}}{\abs{x}^{p(a+1)}}\;
d\mu
\le \int\,\tfrac{\abs{\nabla u}^{p}}{\abs{x}^{pa}}\;d\mu 
+ \Big(\tfrac{\abs{d-p(a+1)}}{p}\Big)^{p-1}\,
\textrm{sgn}(d-p(a+1))\,  \int\, \abs{u}^{p}\,
\tfrac{(x^{t}Ax)^{\frac{p}{2}}}{\abs{x}^{p(a+1)}}\; d\mu
\]
for all $u\in \textrm{C}_{c,0}^{\infty}$. Moreover, if either $A\equiv0$ or $A$ is positive definite,
then the constant $\Big(\tfrac{\abs{d-p(a+1)}}{p}\Big)^{p}$ is optimal (including the
case $p(a+1)=d$). 
\end{thm}

\begin{proofthm}{\ref{caffarelli:1}}
We take $\lambda\ge0$, set $F(x)= \lambda\,\textrm{sgn}(d-p(a+1))\,
\frac{x}{\abs{x}^{p(a+1)}}\,\rho(x)$, and proceed analogously as in the proof of Theorem~\ref{thm-hardy-inequality}.
\end{proofthm}

\begin{thm}\label{caffarelli:2}
Let $A\in \R^{d\times d}$ be a symmetric positive
semi-definite matrix, and let $\beta\in\R$. Then 
\[
\Big(\tfrac{\abs{d-(p+\beta)}}{p}\Big)^{p}\,\int\, 
\tfrac{\abs{u}^{p}}{\abs{x}^{p+\beta}}\;d\mu
\le \int\,\tfrac{\abs{\nabla u}^{p}}{\abs{x}^{\beta}}\;d\mu 
+ \Big(\tfrac{\abs{d-(p+\beta)}}{p}\Big)^{p-1}\,
\textrm{sgn}(d-(p+\beta))\,  \int\, \abs{u}^{p}\,
\tfrac{(x^{t}Ax)^{\frac{p}{2}}}{\abs{x}^{p+\beta}}\; d\mu
\]
for all $u\in \textrm{C}_{c,0}^{\infty}$. Moreover, if either $A\equiv0$ or $A$ is
positive definite, then the constant $\Big(\tfrac{\abs{d-p(a+1)}}{p}\Big)^{p}$ is
optimal (including the case $p(a+1)=d$).
\end{thm}

\begin{proofthm}{\ref{caffarelli:2}}
For $\lambda\ge0$, we set 
$
F(x)= \lambda\,\textrm{sgn}(d-(p+\beta))\, \frac{x}{\abs{x}^{p+\beta}}\,\rho(x),
$
and proceed analogously as in the proof of Theorem~\ref{thm-hardy-inequality}.
\end{proofthm}
%
%
\section{{\bf Application: A nonexistence result}}

\mbox{}\indent In this section, we prove existence and nonexistence of nonnegative
solutions to the p-Kolmogorov parabolic problem~\eqref{IVP}. We define solutions
in the following sense. Similar definitions can be found for example in~\cite{MR2016679} or~\cite{MR1230384}.

\begin{defn}
Let $u_{0}\in L^{2}_{\mu}(]0,+\infty[)$, $f\in L^{2}(0,T;L^{2}_{\mu}(]0,+\infty[))$,
and let $\lambda>0$. We call a function
\[
u \in \textrm{C}([0,T];L^{2}_{\mu}(]0,+\infty[))\cap 
L^{p}(0,T;W^{1,p}_{\mu,loc}(]0,+\infty[))
\]
a \emph{weak solution locally off of zero} of equation
\begin{equation}\label{CP:1-Eq}
\partial_{t}u-\partial_{x}\{\abs{\partial_{x}u}^{p-2}\partial_{x}u\}
=\abs{\partial_{x}u}^{p-2}\partial_{x}u\tfrac{\partial_{x}\rho}{\rho}
+\tfrac{\lambda}{x^{p}}\abs{u}^{p-2}u+f\qquad\text{on $]0,+\infty[\times ]0,T[$,}
\end{equation}
if for all $\mathcal{K}\Subset]0,+\infty[$, $t_{1}, t_{2}\in [0,T] :\,t_{1}\le t_{2}$, 
and all $\varphi\in
W^{1,2}(t_{1},t_{2};L^{2}_{\mu}(\mathcal{K}))\cap
L^{p}(t_{1},t_{2};W^{1,p}_{\mu,0}(\mathcal{K}))$,
\begin{equation}\label{var:form}
\begin{split}
&(u,\varphi)_{L^{2}_{\mu}(\mathcal{K})}\Big\vert_{t_{1}}^{t_{2}}
+ \int_{t_{1}}^{t_{2}}\int_{\mathcal{K}}
\Big\{-u\,\partial_{t}\varphi + \abs{\partial_{x}u}^{p-2}
\partial_{x}u\partial_{x}\varphi-\tfrac{\lambda}{x^{p}}\abs{u}^{p-2}u\,\varphi \Big\}\,
d\mu\,dt\\
&\hspace{3cm} = \int_{t_{1}}^{t_{2}}\int_{\mathcal{K}} f\,\varphi\,d\mu\,dt \mathpunct{.}
\end{split}
\end{equation}
If such a function $u$ is nonnegative a.e. on $]0,+\infty[\times ]0,T[$, then we call
$u$ a nonnegative weak solution locally off of zero of equation~\eqref{CP:1-Eq}.
We call a function $u$ a
\emph{weak solution locally off of zero of initial value problem}
\begin{equation}\label{IVP}
\begin{cases}
\partial_{t}u-\partial_{x}\{\abs{\partial_{x}u}^{p-2}\partial_{x}u\}
 & = \abs{\partial_{x}u}^{p-2}\partial_{x}u\tfrac{\partial_{x}\rho}{\rho}
 + \tfrac{\lambda}{x^{p}}\abs{u}^{p-2}u+f
\qquad\text{on $]0,+\infty[\times ]0,T[$,}\\
\hspace{2,8cm}u(0) & = u_{0}\hspace{4,85cm}\text{on $]0,+\infty[$,}
\end{cases}
\end{equation}
if $u$ is a weak solution locally off of zero of equation~\eqref{CP:1-Eq} and satisfies
$u(0)=u_{0}$ in $L^{2}_{\mu}(]0,+\infty[)$. Here we stress that since
$u\in \textrm{C}([0,+\infty[;L^{2}_{\mu}(]0,+\infty[)$, the initial condition
$u(0)=u_{0}$ in $L^{2}_{\mu}(]0,+\infty[)$ has a meaning.
We call a weak solutions $u$ locally off of zero of initial value problem~\eqref{IVP} \emph{global} if 
\[
u \in \textrm{C}([0,+\infty[;L^{2}_{\mu}(]0,+\infty[))\cap 
L^{p}_{loc}(0,+\infty;W^{1,p}_{\mu,loc}(]0,+\infty[))\mathpunct{.}
\]
\end{defn}

\begin{thm}\label{nonexistence}
Then the following assertions are true:
\begin{enumerate}
\item[(i)] If $0\le\lambda\le C(1,p)$ and if $1<p<+\infty$, then for every nonnegative
$u_{0}\in L^{2}_{\mu}(]0,+\infty[)$ and for every
nonnegative $f\in L^{2}(0,+\infty;L^{2}_{\mu}(]0,+\infty[))$, there is at least one
global nonnegative weak solution $u$ off of zero of initial value problem~\eqref{IVP}
satisfying
\begin{equation}\label{growth-cond-in-d:1}
\norm{u(t)}_{L^{2}_{\mu}(]0,+\infty[)}\le \norm{u_{0}}_{L^{2}_{\mu}}+
\int_{0}^{t}\norm{f(s)}_{L^{2}_{\mu}(]0,+\infty[)}\,ds\qquad\text{for all $t\ge0$.}
\end{equation}
\item[(ii)] If $\lambda>C(1,p)$, if $1<p<2$ and if $f\equiv 0$, then
for any nonnegative $u_{0}\in L^{2}_{\mu}(]0,+\infty[)\setminus\{0\}$, there is no
global nonnegative weak solution $u$ off of zero of initial value problem~\eqref{IVP},
which is bounded with values in $L^{2}_{\mu}(]0,+\infty[)$.
\end{enumerate}
\end{thm}

%
%

\mbox\indent In order to prove the second assertion of
Theorem~\ref{nonexistence}, we need to introduce the Steklov average
of a function $v\in L^{1}(]a,b[\times ]\tau_{1},\tau_{2}[)$ defined on 
$]a,b[\times ]\tau_{1},\tau_{2}[\subseteq \R^{2}$: if
$0<\varepsilon< \tau_{2}-\tau_{1}$ and if $0<h<\varepsilon$, then the Steklov average of $v$ is given by $v_{h}(x,t):=\tfrac{1}{h}\int_{t}^{t+h}v(x,s)\,ds$ for all 
$t\in ]\tau_{1},\tau_{2}-\varepsilon[$, a.e. $x\in ]a,b[$. Moreover, For $q,r\ge 1$
we denote by $L^{q,r}(]a,b[\times ]\tau_{1},\tau_{2}[)$ the parabolic Lebesgue
space $L^{r}(\tau_{1},\tau_{2};L^{q}(]a,b[))$ equipped with the
norm
\[
\norm{u}_{L^{q,r}(]a,b[\times ]\tau_{1},\tau_{2}[)}:=\left(\int_{\tau_{1}}^{\tau_{2}}
\left(\int_{a}^{b}\abs{u(x,t)}^{q}\,dx\right)^{\frac{r}{q}}\, dt\right)^{\frac{1}{r}}\qquad
\text{for all $u\in L^{q,r}(]a,b[\times ]\tau_{1},\tau_{2}[)$.}
\]
Furthermore, we need the following two well-known lemmas. For a proof, we refer the interested reader, for instance, to the book~\cite{MR1230384} of DiBenedetto.

\begin{lem}\label{lemma-of-steklov-averages}
Let $]a,b[, ]\tau_{1},\tau_{2}[\subseteq\R$ be two open intervals. Then the following assertions are valid.
\begin{enumerate}
\item[(i)] If $v\in L^{q,r}(]a,b[\times ]\tau_{1},\tau_{2}[)$, then for every 
$0<\varepsilon<\tau_{2}-\tau_{1}$,
\[
v_{h}\to v\qquad\text{in $L^{q,r}(]a,b[\times ]\tau_{1},\tau_{2}-\varepsilon[)$}
\quad\text{as $h\to0+$.}
\]
\item[(ii)] If $v\in \textrm{C}([\tau_{1},\tau_{2}];L^{q}(]a,b[))$, then $v_{h}(t)$ can 
be defined in $t=\tau_{1}$ by $v_{h}(x,\tau_{1})=\tfrac{1}{h}
\int_{\tau_{1}}^{\tau_{1}+h}v(x,s)\,ds$ $(x\in ]a,b[)$ for all $0<h<\tau_{2}-\tau_{1}$. 
Moreover, then for every $0<\varepsilon<\tau_{2}-\tau_{1}$ and every 
$t\in [\tau_{1},\tau_{2}-\varepsilon[$,
\[
v_{h}(t)\to v(t)\qquad\text{in $L^{q}(]a,b[)$}\quad\text{as $h\to0+$.}
\]
\end{enumerate}
\end{lem}

\begin{lem}\label{lem:integration-by-parts}
If $u$ is a weak solution locally off of zero of equation~\eqref{CP:1-Eq} and if 
$g : \R\to \R$ is Lipschitz-continuous, then for every
$\phi\in W^{1,p}_{\mu}(]0,+\infty[)\cap L^{2}_{\mu}(]0,+\infty[)$ having compact
support contained in $]0,+\infty[$, and every $t_{1}, t_{2}\in [0,T]:\,t_{1}<t_{2}$,
\begin{equation}\label{int:g-equality}
\begin{split}
&\int_{0}^{+\infty}\left(\int_{0}^{u(t_{2})}g(s)\,ds\right)\,\phi\, d\mu
-\int_{0}^{+\infty}\left(\int_{0}^{u(t_{1})}g(s)\,ds\right)\,\phi\, d\mu\\
&\qquad + \int_{t_{1}}^{t_{2}}\int_{0}^{+\infty} \abs{\partial_{x}u}^{p}\, g'(u)\, 
\phi\,d\mu\,dt
+\int_{t_{1}}^{t_{2}}\int_{0}^{+\infty} \abs{\partial_{x}u}^{p-2}\partial_{x}u\, g(u)\,
 \partial_{x}\phi\,d\mu\,dt\\
&\qquad\qquad = \int_{t_{1}}^{t_{2}}\int_{0}^{+\infty} 
\left(\tfrac{\lambda}{x^{p}}\abs{u}^{p-2}u+f\right)\,g(u)\, \phi\,d\mu\,dt \mathpunct{.}
\end{split}
\end{equation}
\end{lem}

\begin{myproof}
By a standard approximation argument, one sees that it is sufficient to prove
the claim of this lemma for test-functions $\phi\in \textrm{C}_{c}^{1}(]0,+\infty[)$.
Thus, we fix $\phi\in \textrm{C}_{c}^{1}(]0,+\infty[)$, and for fixed $0<t<t+h<T$,
we take $t_{1}=t$, $t_{2}=t+h$, and multiply equation~\eqref{var:form} by 
$h^{-1}$. Then equation~\eqref{var:form} becomes
\[
\begin{split}
&\int_{0}^{+\infty}h^{-1}\left( u(t+h)-u(t)\right)\,\phi\,d\mu 
+ \int_{t}^{t+h}\int_{0}^{+\infty}\,h^{-1}\left(\abs{\partial_{x}u}^{p-2}
\partial_{x}u\right)\, \partial_{x}\phi\,d\mu\,dt\\
&\qquad\qquad= \int_{t}^{t+h}\int_{0}^{+\infty}\, h^{-1}
\left(\tfrac{\lambda}{x^{p}}\,\abs{u}^{p-2}u+f\right)\,\phi \,d\mu\,dt\mathpunct{.}
\end{split}
\]
By Fubini's theorem, since $\partial_{t}u_{h}(t)=h^{-1}\left( u(t+h)-u(t)\right)$, and
by the definition of Steklov averages, the last equality can be rewritten as
\[
\int_{0}^{+\infty}\partial_{t}u_{h}(t)\,\phi\,d\mu 
+ \int_{0}^{+\infty}\,\left(\abs{\partial_{x}u}^{p-2}
\partial_{x}u\right)_{h}(t)\, \partial_{x}\phi\,d\mu
= \int_{0}^{+\infty}\, 
\left(\tfrac{\lambda}{x^{p}}\abs{u}^{p-2}u+f\right)_{h}(t)\,\phi \,d\mu\mathpunct{.}
\]
By the hypothesis, for any $t\in ]0,T[$, $g(u_{h}(t))\,\phi\in 
W^{1,p}_{\mu,0}(\mathcal{K})$ with distributional derivative
\[
\partial_{x}\left(g(u_{h}(t))\,\phi \right)= g'(u_{h}(t))(\partial_{x}u)_{h}(t)
\,\phi+g(u_{h}(t))\,\partial_{x}\phi\mathpunct{,}
\]
where $\mathcal{K}\Subset]0,+\infty[$ can be chosen as an open and
bounded interval such that the support of $\phi$
$\textrm{supp}(\phi)\subseteq \mathcal{K}$. Thus, we can replace 
$\phi$ by $g(u_{h}(t))\,\phi$ in the last
equality. Then, integrating over $]t_{1},t_{2}[$ for any
$t_{1}, t_{2}\in [0,T]:\,t_{1}<t_{2}$ and applying Fubini's theorem gives
\[
\begin{split}
&\int_{0}^{+\infty}\left(\int_{0}^{u_{h}(t_{2})}g(s)\,ds\right)\,\phi\,d\mu-  
\int_{0}^{+\infty}\left(\int_{0}^{u_{h}(t_{1})}g(s)\,ds\right)\,\phi\,d\mu\\
& + \int_{t_{1}}^{t_{2}}\int_{0}^{+\infty}\,\left\{\left(\abs{\partial_{x}u}^{p-2}
\partial_{x}u\right)_{h}(t)\, g'(u_{h}(t))(\partial_{x}u)_{h}(t)
\,\phi\,+ \left(\abs{\partial_{x}u}^{p-2}\partial_{x}u\right)_{h}(t)\,g(u_{h}(t))\,
\partial_{x}\phi\right\} d\mu\,dt\\
&\hspace{2cm} = \int_{t_{1}}^{t_{2}}\int_{0}^{+\infty}\, 
\left(\tfrac{\lambda}{x^{p}}\abs{u}^{p-2}u+f\right)_{h}(t)\,g(u_{h}(t))\,\phi\,d\mu\,dt
\mathpunct{.}
\end{split}
\]
Now, by sending $h\to 0+$ in the last equality and using
Lemma~\ref{lemma-of-steklov-averages} leads to equality~\eqref{int:g-equality}. 
\end{myproof}

\begin{proofthm}{\ref{nonexistence}}
If $\lambda\le C(1,p)$, then by Theorem~3.17 in
\cite{goldstein-hauer-rhandi-2012:1}, for every nonnegative
$u_{0}\in L^{2}_{\mu}(]0,+\infty[)$ and for every
nonnegative $f\in L^{2}(0,+\infty;L^{2}_{\mu}(]0,+\infty[))$, there exists at least one
global nonnegative weak solution $u$ off of zero of initial value problem~\eqref{IVP}.
Moreover, there is a sequence $(u_{m})_{m\ge1}$ of strong solutions 
\[
u_{m}\in W^{1,2}_{loc}(]0,+\infty[;L^{2}_{\mu}(]0,+\infty[))\cap 
L^{p}_{loc}(0,+\infty;W^{1,p}_{\mu}(]0,+\infty[))
\]
of the truncated problem
\begin{align}\label{CP:6-Eq}\tag{$1_{m}$}
&\; \partial_{t}u_{m}-\partial_{x}\{\abs{\partial_{x}u_{m}}^{p-2}\partial_{x}u_{m}\}
 = \abs{\partial_{x}u_{m}}^{p-2}\partial_{x}u_{m}\tfrac{\partial_{x}\rho}{\rho}
+ \Phi_{m}\,u^{p-1}_{m}+f_{m}
&&\text{on $]0,+\infty[\times ]0,+\infty[$,}\\
\label{CP:6-IC}\tag{$2_{m}$}
&\; u_{m}(0) =u_{m}^{0} &&\text{on $]0,+\infty[$,}
\end{align}
for $\Phi_{m}(x):=\min\{\tfrac{\lambda}{x^{p}},m\}$, $f_{m}(x,t)=\min\{f(x,t),m\}$,
and $u_{m}^{0}(x):=\min\{u_{0}(x),m\}$, each of them defined pointwise almost everywhere, such that for every $T>0$,
\begin{equation}\label{strong-limit:um}
u_{m}\to u\qquad\text{in $\textrm{C}([0,T];L^{2}_{\mu}(]0,+\infty[))$
as $m\to+\infty$.}
\end{equation}
If we first multiply equation~\eqref{CP:6-Eq} by $u_{m}$ with respect to the inner
product on $L^{2}_{\mu}(]0,+\infty[)$, and subsequently integrate from $0$ to $t$
for fixed $t\ge0$, then we obtain by Hardy's inequality~\eqref{hardy-inequality:1}
and by the Cauchy-Schwarz inequality that
\[
\tfrac{1}{2}\norm{u_{m}(t)}_{L^{2}_{\mu}(]0,+\infty[)}^{2}
\le \tfrac{1}{2}\norm{u_{m}^{0}}_{L^{2}_{\mu}(]0,+\infty[)}^{2}+
\int_{0}^{t}\norm{f_{m}(s)}_{L^{2}_{\mu}(]0,+\infty[)}\,
\norm{u_{m}(s)}_{L^{2}_{\mu}(]0,+\infty[)}\,ds\mathpunct{.}
\]
By Lemma A.5 in the book~\cite{MR0348562} by Brezis, we can deduce from the last inequality that
\[
\norm{u_{m}(t)}_{L^{2}_{\mu}(]0,+\infty[)}\le \norm{u_{0}}_{L^{2}_{\mu}}+
\int_{0}^{t}\norm{f(s)}_{L^{2}_{\mu}(]0,+\infty[)}\,ds\qquad\text{for all $t\ge0$.}
\]
Since  $\norm{u_{m}^{0}}_{L^{2}_{\mu}(]0,+\infty[)}\le
\norm{u_{0}}_{L^{2}_{\mu}(]0,+\infty[)}$ and
$\norm{f_{m}}_{L^{2}(0,T,L^{2}_{\mu}(]0,+\infty[))}\le
\norm{f}_{L^{2}(0,T,L^{2}_{\mu}(]0,+\infty[))}$, and by
limit~\eqref{strong-limit:um}, we can send $m\to+\infty$ in the last inequality and
see that $u$ satisfies inequality~\eqref{growth-cond-in-d:1}. Thus claim~(i) 
of this theorem holds true.\\[7pt]
%
%
\mbox{}\indent Let $1<p<2$ and $f\equiv0$. We suppose that there is a
global nonnegative weak solution $u$ off of zero of equation~\eqref{CP:1-Eq}, 
with nonnegative  initial value $u(0)=u_{0}\in L^{2}_{\mu}(]0,+\infty[)\setminus\{0\}$,
which is bounded with values in $L^{2}_{\mu}(]0,+\infty[)$. Then, we shall
reach a contradiction. To see this, we fix 
$\varphi\in \textrm{C}_{c}^{\infty}(]0,+\infty[)$ and for every integer $k\ge1$ and
every $s\in \R$, let $g_{k}(s)=(s+1/k)^{1-p}$. Then, by
Lemma~\ref{lem:integration-by-parts} for $g=g_{k}$,
$\phi=\abs{\varphi}^{p}$, $t_{1}=0$, and $t_{2}=t$ for any fixed $t\ge0$, 
\[
\begin{split}
&\tfrac{1}{2-p}\int_{0}^{+\infty}(u(t)+1/k)^{2-p}\,\abs{\varphi}^{p}\, d\mu
-\tfrac{1}{2-p}\int_{0}^{+\infty}(u_{0}+1/k)^{2-p}\,\abs{\varphi}^{p}\, d\mu\\
&\hspace{2cm} + (1-p)\int_{0}^{t}\int_{0}^{+\infty} \abs{\partial_{x}u(s)}^{p}\,
 (u(s)+1/k)^{-p}\, \abs{\varphi}^{p}\,d\mu\,ds\\
&\hspace{2cm} + p\,\int_{0}^{t}\int_{0}^{+\infty} \abs{\partial_{x}u(s)}^{p-2}
\partial_{x}u(s)\,(u(s)+1/k)^{1-p}\,\partial_{x}\varphi\,\abs{\varphi}^{p-2}
\varphi\,d\mu\,ds\\
&\hspace{4cm} = \int_{0}^{t}\int_{0}^{+\infty} 
\tfrac{\lambda}{x^{p}}u^{p-1}(s)\,(u(s)+1/k)^{1-p}\, \abs{\varphi}^{p}\,d\mu\,ds 
\mathpunct{.}
\end{split}
\]
We apply Young's inequality, and since $(u_{0}+1/k)^{2-p}\,\abs{\varphi}^{p}$ is nonnegative a.e. on $]0,+\infty[$, we obtain that 
\begin{equation}\label{previous-inequality}
\begin{split}
&\int_{0}^{t}\int_{0}^{+\infty} 
\tfrac{\lambda}{x^{p}}u^{p-1}(s)\,(u(s)+1/k)^{1-p}\, \abs{\varphi}^{p}\,d\mu\,ds\\
&\hspace{2cm} \le t\,\int_{0}^{+\infty}\abs{\partial_{x}\varphi}^{p}\,d\mu
 +  \tfrac{1}{2-p}\int_{0}^{+\infty}(u(t)+1/k)^{2-p}\,\abs{\varphi}^{p}\, d\mu\mathpunct{.}
\end{split}
\end{equation}
For almost every $(x,s)\in ]0,+\infty[\times ]0,t[$,
\[
0\le \tfrac{\lambda}{\abs{x}^{p}}\,\tfrac{u^{p-1}(x,s)\,\abs{\varphi(x)}^{p}}{(u(x,s)+\frac{1}{k})^{p-1}}\nearrow \tfrac{\lambda}{\abs{x}^{p}}\,\abs{\varphi(x)}^{p}\qquad\text{as $k\to+\infty$,}
\]
and
\[
(u(x,s)+\tfrac{1}{k})^{2-p}\,\abs{\varphi(x)}^{p} \searrow u^{2-p}(x,s)\,\abs{\varphi(x)}^{p}\qquad\text{as $k\to+\infty$.}
\]
Thus, by Beppo-Levi's convergence theorem, sending $k\to+\infty$ in
inequality~\eqref{previous-inequality} gives
\[
t\,\int_{0}^{+\infty} \tfrac{\lambda}{x^{p}}\, \abs{\varphi}^{p}\,d\mu
- t\,\int_{0}^{+\infty}\abs{\partial_{x}\varphi}^{p}\,d\mu
\le \tfrac{1}{2-p}\int_{0}^{+\infty}u^{2-p}(t)\,\abs{\varphi}^{p}\, d\mu\mathpunct{,}
\]
and by H\"older's inequality, 
\[
t\,\int_{0}^{+\infty} \tfrac{\lambda}{x^{p}}\, \abs{\varphi}^{p}\,d\mu\,
- t\,\int_{0}^{+\infty}\abs{\partial_{x}\varphi}^{p}\,d\mu
\le \tfrac{1}{2-p}
\le \norm{\varphi}_{L^{2}(]0,+\infty[)}^{p}\,\norm{u(t)}_{L^{2}_{\mu}(]0,+\infty[)}^{2-p}
\mathpunct{.}
\]
We divide this inequality by $t>0$. Since $t\mapsto u(t)$ is bounded from
$]0,+\infty[$ to $L^{2}_{\mu}(]0,+\infty[)$, sending $t+\infty$ yields to
\[
\int_{0}^{+\infty} \tfrac{\lambda}{x^{p}}\, \abs{\varphi}^{p}\,d\mu\,
- \int_{0}^{+\infty}\abs{\partial_{x}\varphi}^{p}\,d\mu
\le 0\mathpunct{.}
\]
Since $\varphi\in \textrm{C}_{c}^{\infty}(]0,+\infty[)$ has been arbitrary in this
inequality and since $\textrm{C}_{c}^{\infty}(]0,+\infty[)$ lies dense in
$W^{1,p}_{\mu}(]0,+\infty[) $, we have thereby shown that
\[
\inf_{\varphi \in W^{1,p}_{\mu}(]0,+\infty[): 
\norm{\varphi}_{L^{p}_{\mu}(]0,+\infty[)}>0}
\frac{\displaystyle\int_{0}^{+\infty}\,\abs{\partial_{x}\varphi}^{p}\,d\mu 
- \displaystyle\int_{0}^{+\infty}\, \frac{\lambda}{x^{p}}\abs{\varphi}^{p} \, d\mu}{\displaystyle\int_{0}^{+\infty}\,\abs{\varphi}^{p}\,d\mu}\ge0\mathpunct{,}
\]
but this obviously contradicts the optimality of the constant $C(1,p)$.
\end{proofthm}
%
%
%

\providecommand{\bysame}{\leavevmode\hbox to3em{\hrulefill}\thinspace}
\providecommand{\MR}{\relax\ifhmode\unskip\space\fi MR }
\providecommand{\MRhref}[2]{%
  \href{http://www.ams.org/mathscinet-getitem?mr=#1}{#2}
}
\providecommand{\href}[2]{#2}

 

\mbox{}
\vfill

\end{document}